\documentclass[twoside,final,reqno]{amsart}

\usepackage{amsfonts}
\usepackage{amssymb}
\usepackage{latexsym}
\usepackage{verbatim}
\usepackage{epsfig}

\newtheorem{theorem}{Theorem}[section]

\newtheorem{corollary}[theorem]{Corollary}
\newtheorem{proposition}[theorem]{Proposition}

\theoremstyle{definition}
\newtheorem{remark}[theorem]{Remark}
\newtheorem{remarks}[theorem]{Remarks}
\newtheorem{example}[theorem]{Example}
\newtheorem{definition}[theorem]{Definition}

\newtheorem{note}[theorem]{Note}

\newenvironment{Proof}{\removelastskip\par\medskip
\noindent{\em Proof.} \rm}{\penalty-20\null\hfill$\square$\par\medbreak}

\renewcommand{\r}{{\rm{rhd\hspace{2pt}}}}
\newcommand{\rHd}{{\rm{rHd\hspace{2pt}}}}
\newcommand{\h}{{\rm{ht}}}

\newcommand{\inter}{{\rm{int}}}

\newcommand{\Sing}{{\rm{Sing\hspace{2pt}}}}

\newcommand{\id}{{\rm{id}}}
\newcommand{\Jac}{{\rm{Jac}}}

\newcommand{\im}{{\rm{Im}}}

\newcommand{\Iso}{{\rm{Iso}}}
\newcommand{\cl}{{\rm{closure}}}

\newcommand{\ity}{{\infty}}

\newcommand{\e}{\varepsilon}
\newcommand{\fin}{\hspace*{\fill}$\square$\vspace*{2mm}}

\newcommand{\cA}{{\mathcal A}}

\newcommand{\cG}{{\mathcal G}}

\newcommand{\cW}{{\mathcal W}}

\newcommand{\cN}{{\mathcal N}}
\newcommand{\cY}{{\mathcal Y}}

\newcommand{\bC}{{\mathbb C}}

\newcommand{\bP}{{\mathbb P}}

\newcommand{\bX}{{\mathbb X}}

\newcommand{\bY}{{\mathbb Y}}

\begin{document}

\title[Singularities and Topology of Meromorphic Functions]
 {Singularities and Topology 
 of Meromorphic \\ Functions}

\author[Mihai Tib\u ar]{Mihai Tib\u ar}

\address{Math\' ematiques, UMR 8524 CNRS, 
Universit\'e des Sciences et Tech. de Lille,  
 59655 Villeneuve d'Ascq, France.}

\email{tibar@agat.univ-lille1.fr}

\thanks{Partially supported by the Newton Institute at Cambridge and by the European Commission and the Institute of Mathematics of the Romanian Academy under the EURROMMAT contract ICA1-CT-2000-70022.
}



\begin{abstract}
We present several aspects of the ``topology of meromorphic functions'', which we conceive as a general theory which includes the topology of holomorphic functions, the topology of pencils on quasi-projective spaces and the topology of polynomial functions. 
\end{abstract}

\maketitle

\tableofcontents
\setcounter{section}{0}
\section{Introduction}

Milnor \cite{Mi} defined the basic ingredients for studying the topology of holomorphic germs of functions $f \colon (\bC^n,0)\to \bC$. 
One of the main motivations for Milnor's book concerns the link $K$ of an isolated singularity and its complement $S^{2n-1} \setminus K$: links which are exotic speres have been discovered by Hirzebruch \cite{Hi} and Brieskorn \cite{Bri}; in case $n=2$, the components of $K$ are iterated toric knots. The study of isolated singularities ever since revealed striking phenomena and established bridges between several branches of mathematics.

In another stream of research, there has been an increasing interest in the last decade
 for the study of the global topology of polynomial functions, especially
 in connection with the behaviour at infinity. This topic is closely related to the affine geometry and to dynamical systems on non-compact spaces.

This paper reports on how to extend the study from holomorphic germs to the class of {\em meromorphic functions}, local or global.  In the same time, polynomial functions $\bC^n \to \bC$ can be viewed as a special case of global meromorphic functions, as we explain futher on. A global meromorphic function defines a pencil of hypersurfaces, therefore our approach moreover yields a generalization of the theory of {\em Lefschetz pencils}.

 Another motivation for the study of meromorphic functions is Arnold's approach to the classification of simple germs of 
meromorphic functions under certain equivalence relations \cite{Ar}.

\begin{figure}[hbtp]\label{f:spec}
\begin{center}
\epsfxsize=7cm
\leavevmode
\epsffile{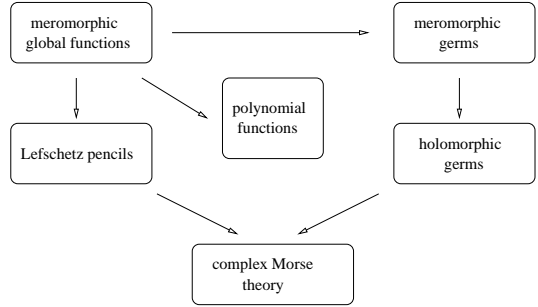}
\end{center}
\caption{{\em 
 Specialization of topics}}
\label{f:1}   
\end{figure}

 Let us introduce the first definitions.
A meromorphic function, or pencil of hypersurfaces, on a compact complex analytic space $Y$, is a function $F \colon Y \dashrightarrow \bP^1$ defined as the ratio of two sections $P$ 
and $Q$ of a holomorphic line bundle over $Y$.  Then $F = P/Q$ is 
 a holomorphic function on $Y\setminus A$, where $A := \{ P=Q=0\}$ is the {\em  base locus} of the pencil (also called {\em  axis}, or {\em indeterminacy locus}).
 A germ of meromorphic function on a space germ is just the ratio of two holomorphic germs $f = p/q  \colon (\cY, y) \dashrightarrow \bP^1$. By definition, $f$ is equal to $f'=p'/q'$, as germs at $y$, if and only if there exists a holomorphic germ $u$ such that $u(y)\not= 0$ and that $p=up'$, $q=uq'$. Then $f$ is holomorphic on the germ at $y$ of the complement $\cY \setminus \cA$ of the axis $\cA = \{ p=q=0\}$.

  Meromorphic functions give rise to a new type of singularities, those occuring along the indeterminacy locus. To define them, we need to introduce some more objects attached to a meromorphic function.
  
\begin{definition}
Let $G := \{ (x, \tau)\in (Y\setminus A) \times \bP^1 \mid F(x) =\tau\}$
and let $\bY$ denote the analytic closure of $G$ in  $Y\times \bP^1$, namely:
\[ \bY =  \{ (x,[s:t])\in Y\times \bP^1 \mid tP(x) - sQ(x) =0\}.\]
  In case of a germ of meromorphic function at $(\cY,y)$, one similarly  defines space germs, which we denote by $(G,(y,\tau))$, resp. $(\bY,(y,\tau))$. See also Definition \ref{d:local}.
\end{definition}
  
 First note that $G$ is the graph of the restriction $F_{| Y\setminus A}$. Therefore 
 $Y\setminus A \simeq G$ embeds into $\bY$ and the projection
 $\pi : \bY \to \bP^1$ is an extension of the function $F_{| Y\setminus A}$.
  One may also say that $\bY \stackrel{\sigma}{\rightarrow} Y$ is a blow-up of $Y$ along the axis $A$, such that the meromorphic function $F \colon Y \dashrightarrow \bP^1$ pulls back to a well defined holomorphic function $\pi : \bY \to \bP^1$.
   
\begin{equation}\label{eq:blow}
 \begin{array}{ccc}
\bY & \  & \  \\
 \mbox{\tiny $\sigma$} \downarrow \ \ &  \ \  \ \  \searrow \mbox{\tiny $\pi$}& \\
 Y &  \stackrel{\mbox{\tiny $F$}}{\dashrightarrow} & \bP^1
\end{array}  
\end{equation}  

  We shall also consider the restriction of  $F$ (or of a germ $f$), to $X:= Y\setminus V$, where $V$ is some compact analytic subspace of $Y$. The case $V= \{ Q=0\}$ is of particular interest for the following reason.
  Let $P \colon \bC^n \to \bC$ be a polynomial of degree $d$, let $\tilde P$ be the homogenized of $P$ by the new variable $x_0$ and let $H^\ity =\{ x_0 =0\}$ be the hyperplane at infinity. Then 
  $\tilde P/x_0^d \colon \bP^n \dashrightarrow \bP^1$ is a meromorphic function on $Y := \bP^n$ which coincides with $P$ over $\bP^n\setminus H^\ity$. We shall briefly outline in \S \ref{s:poly} some results and literature on polynomials.

 Our meromorphic function (as a global one or as a germ) defines a family (=pencil) of hypersurfaces
  on each of the spaces defined above:  $\bY$, $Y$, $Y\setminus A$,  or $X=Y\setminus V$. The map to $\bP^1$, whenever defined, yields the pencil, as the family of its fibres. In other cases, we take the closures of fibres in the considered space. For instance, in case of $Y$, we take the closure
  of each hypersurface $F_{|Y\setminus A}^{-1}(\tau)$ within $Y$, for $\tau\in \bP^1$;
 each such closure contains $A$ and we actually have  $\overline{F_{|Y\setminus A}^{-1}(\tau)}= \pi^{-1}(\tau)$.
 The role of the ``completed space" $\bY$ is that it contains all these  pencils: we just restrict the fibres of $\pi$ to the particular  
 subspace of $Y$ and we get the pencil we are looking for.
 
 With this approach, one covers a large field. For instance, the class of holomorphic functions (or germs) represents the case when $A=\emptyset$.

 Now, for any pencil on $Y$ or on $X$, there is only a finite number of atypical values, or atypical fibres. This finiteness result comes from ideas of Thom \cite{Th} and is based on the fact that we can stratify the space $Y$ (such that $V$ is union of strata, in case $V\not= \emptyset$), restrict to $Y\setminus A \simeq \bY\setminus (A\times \bP^1)$ and extend this to some Whitney stratification of $\bY$. 

 
 In case of a germ at $(\cY,y)$ of a meromorphic function, one considers the germ of such a Whitney stratification at $\{y\}\times \bP^1 \subset \bY$. Local finiteness of the strata implies that non-transversality of the projection $\pi$ happens at discrete values only.
 
\begin{proposition}\label{p:1} 
There exists a finite set $\Lambda\subset\bP^1$ such that the map
$\pi : \bY \setminus \pi^{-1}(\Lambda) \to 
\bP^1\setminus \Lambda$ is a stratified
locally trivial {\rm C}$^0$ fibration.

In particular, the restrictions $\pi_| \colon \bY \setminus ((V\times \bP^1) \cup \pi^{-1}(\Lambda)) \to \bP^1 \setminus \Lambda$ and $F_|\colon
Y \setminus (V\cap A\cap F^{-1}(\Lambda)) \to \bP^1\setminus \Lambda$  
 are  stratified
locally trivial fibrations. 
\fin 
\end{proposition} 

In our approach, the singularities of the meromorphic function $F$ along the indeterminacy locus $A$  
are the stratified singularities of $\pi$ at $(A\times \bP^1)\cap \bY$. Usually, singularities of functions on singular spaces
are defined with respect to some Whitney stratification. Here we use
 instead a {\em partial Thom stratification}, denoted $\cG$, as we already used in particular cases (cf. \cite{Ti-t,Ti-lef}, \cite{ST-x}). This is a more general type of stratification
since Whitney (b) condition is not required. Nevertheless, it allows one to study topological aspects, including homotopy type, at least for isolated singularities, in both local or global context.

Instead of endowing $\bY$ with a stratified structure, another strategy for studying the topology of the meromorphic function $F$
would be to further blow up $\bY$ in diagram (\ref{eq:blow}), such that the pull-back of $\{ P=0\} \cup \{ Q=0\}$ becomes a divisor with normal crossings. One may then use the data provided by this divisor in order to 
get informations. In this spirit, some results were found in the polynomial case, in two variables, by Fourrier \cite{Fo} and L\^e-Weber \cite{LW};  computation of the zeta functions of the monodromy has been done for polynomials and particular meromorphic germs (namely for $Y$ nonsingular and $X = Y\setminus A$) by Gusein-Zade, Melle and Luengo \cite{GLM, GLM-3, GLM-4}.

This paper revisits the techniques and results of \cite{Ti-t,Ti-lef}, \cite{ST-x} and extends them to the more general context introduced above. The main scope is to show how to study  
vanishing cycles of meromorphic functions in both 
local and global context. 



\section{Vanishing homology and singularities}\label{s:van}


Let us first define vanishing homology attached to a global meromorphic function and relate it to the singularities along the indeterminacy locus.
The vanishing homology is important in detecting and controlling (whenever possible) the change of topology of the fibres.
  
We shall use the following notations. For any subset $W\subset \bP^1$,  $\bY_W := \pi^{-1}(W)$, 
$Y_W :=  A \cup F^{-1}(W)$, $X_W := X\cap Y_W$.
 The special case $X= Y\setminus A$ and this is the object of study, explicitly, in \cite{ST-x}, and implicitly, in \cite{Ti-t}. Our presentation follows the one of \cite{ST-x}, adapting it to our more general situation; in particular our notations are different 
 from those in \cite{ST-x}.
 
 Let $a_i\in \Lambda$ be an atypical value of $\pi$ and take a small enough disc $D_i$ at $a_i$ such that $D_i\cap \Lambda = \{ a_i \}$. 
 Let's fix some point $s_i\in\partial D_i$. Let $s \in \bP^1\setminus \Lambda$ be a general value, situated on the boundary of some  big closed disc $D \subset \bP^1$, such that $D \supset D_i$, $\forall a_i\in D \cap \Lambda$, and that $D \cap D_i =\emptyset$, $\forall a_i\in D \setminus \Lambda$. 
 
 The vanishing homology of meromorphic functions should be a natural
extension of the vanishing homology of local holomorphic functions. In the latter case, the total space of the Milnor fibration \cite{Mi} is contractible, by the local conical structure of analytic sets \cite{BV}. For global meromorphic functions, the total space one may take cannot be contractible anymore and the general fibre $X_s$ inherits 
part of its homology. 

\begin{definition}\label{d:vanhom} 
The  {\em vanishing homology of $F_{|X}$ at $a_i$} is 
the relative homology \\
$H_*(X_{D_i}, X_{s_i})$.
\end{definition}
 In the case $X = Y\setminus A$, this corresponds to the definition 
 used by Siersma and the author in \cite{ST-x}.
 
We identify $X_s$ to $X_{s_i}$, in the following explicit manner.  
For each $i$, take a path $\gamma_i \subset D$ from $s$ 
to $s_i$, with the usual conditions: the 
path $\gamma_i$ has no self intersections and does not intersect any 
other path $\gamma_j$, except at the point $s$.
 Then Proposition \ref{p:1} allows identifying $X_s$ to $X_{s_i}$, 
 by parallel transport along $\gamma_i$.

 A general result tells that vanishing homologies can be ``patched" together.
 This type of result was observed before in different particular situations, see e.g. \cite[\S 5]{Br}, \cite{Si}. More precisely, we have the following result, extending the context of \cite[Proposition 2.1]{ST-x} to any $X= Y\setminus V$:
\begin{proposition}\label{p:basic} \
\begin{enumerate}
\item  $ H_*(X_{D}, X_s) = \oplus_{a_i\in\Lambda} H_*(X_{D_i}, X_{s_i}).$
\item The long exact sequence of the triple $(X_{D}, X_{D_i}, X_s)$
decomposes into short exact sequences which split:
\begin{equation}\label{eq:exact} 
0\to H_*(X_{D_i}, X_{s_i}) \to H_*(X_{D},X_s) \to 
H_*(X_{D}, X_{D_i}) \to 0 .
\end{equation}
 \item There is a natural identification $H_*(X_{D},X_{D_i}) = 
\oplus_{a_j\in\Lambda, j\not= i} H_*(X_{D_j}, X_s)$. 
\end{enumerate}
\end{proposition}
\begin{Proof}
  By Proposition \ref{p:1}, 
the fibration $F_|: \bY \setminus ((v\times \bP^1) \cup \pi^{-1}(\Lambda) \to \bP^1 \setminus 
\Lambda$ is locally trivial. Its fibre over some $b$ is, by definition, $X_b$. We then get a sequence of excisions:
\[\oplus_{a_i\in\Lambda} H_*(X_{D_i}, X_{s_i})  \stackrel{\simeq}{\longrightarrow}
H_*( \pi^{-1}(\cup_{a_i\in \Lambda} D_i\cup\gamma_i), X_s) 
\stackrel{\simeq}{\longrightarrow} H_*(X_{D}, X_s).
\]
This also shows that each inclusion $(X_{D_i}, X_{s_i}) \subset 
(X_{D},X_{s_i})$ 
induces an injection in homology $H_*(X_{D_i}, X_{s_i}) \hookrightarrow H_*(X_{D}, X_s)$.
 The points (a), (b), (c) all follow from this. 
\end{Proof}

Vanishing homology of $F$ has its local counterpart and is closely related to singularities of $F$. One would like to say that vanishing homology is supported at the singular points of $F$. The typical problem for a meromorphic function is that it has singularities also outside the ground space, $Y$ or $X$. Therefore we need a larger space, such as $\bY$, to define singularities. Then the support of vanishing cycles is included into the singular locus of $\pi$ on $\bY$. (In cohomology, the sheaf of vanishing cycles of a function $h$ on a nonsingular space is indeed supported by the singular locus of $h$, see \cite{De}.) 

Let now give the precise definition of what we consider as singularities of $F$ (resp. of $f$). We relax the stratification conditions at $A\times \bP^1 \subset \bY$ and use only the Thom condition. Let us first recall the latter, following \cite{Ti-t}.

Let $\cG = \{ \cG_\alpha\}_{\alpha \in S}$ be a locally finite 
stratification such that $\bY \setminus X$ is union of strata.
Let $\xi :=(y,a)$ be a point on a stratum $\cG_\alpha$.
We assume, without loss of generality, that $a\not= [1:0]$. 
Let $f =p/q$ be our meromorphic germ on $(\cY, y)$ or a local representative of the germ of $F$ at $(Y, y)$. Then $q=0$ is a local equation for $A\times \bP^1$ at $\xi$. The {\em 
Thom regularity condition} (a$_q$) at  $\xi \in \cG_\alpha$
is satisfied (see e.g. 
\cite[ch. I]{GWPL}) if for any stratum $\cG_\beta$ such that $\cG_\alpha \subset \bar \cG_\beta$, the relative conormal space (see \cite{Te-2}, \cite{HMS} for a 
definition) of $q$ on $\bar \cG_\beta$
is included into the conormal of $\cG_\alpha$, locally at $\xi$, i.e.,  $(T_{\cG_\alpha}^*)_\xi \supset (T^*_{q| \overline{ \cG_\beta}})_\xi$.
This condition is known to be independent on $q$, up to 
multiplication by 
a unit \cite[Prop. 3.2]{Ti-t}. 


\begin{definition}\label{d:thom} 
Let $\cG$ be a stratification on $\bY$ as above such that it restricts to a Whitney stratification on $Y\setminus A$, where $V\setminus A$ is union of strata. We say that $\cG$ is a {\em partial 
Thom stratification} ($\partial \tau$-stratification) if Thom's condition (a$_q$) is satisfied at any point $\xi\in A\times \bP^1 \subset \bY$.
\end{definition}

One may extend a Whitney stratification on $Y\setminus A$ to a locally finite $\partial \tau$-stratification of $\bY$, by usual stratification theory arguments (see e.g. 
\cite{GWPL}). For instance, the Whitney stratification $\cW$ of $\bX$ that we have considered before
is an example of $\partial 
\tau$-stratification. This follows from  \cite[Th\'eor\`eme 
4.2.1]{BMM} or \cite[Theorem 3.9]{Ti-t}.
One can also construct a canonical (minimal) $\partial \tau$-stratification;
we send to \cite{Ti-t} for further details.
\hyphenation{stra-ti-fi-cation}
  

\begin{definition}\label{d:sing}
Let $\cG$ be a $\partial \tau$-stratification on $\bY$. We say that the following closed subset of $\bY$:
\[ \Sing_{\cG}F := \cup_{
 \cG_\alpha \in \cG }
 \cl(\Sing \pi_{|\cG_\alpha} )\]
 is the singular locus of $F$ with respect to $\cG$.
We say that $F$ has {\em isolated singularities} with respect to $\cG$ 
if $\dim \Sing_\cG F \le 0$. 
\end{definition}

For the singular locus of a germ $f$, one modifies this definition accordingly.
  The singularities of the new type are those along the indeterminacy locus, namely on  $A\times \bP^1$. We shall further investigate the relation between singularities and 
  vanishing homology, in case of isolated singularities.
Let us define the local fibration of a meromorphic germ, already used in particular cases in \cite{ST}, \cite{Ti-t}, \cite{ST-x}), and which relates to the fibration of a holomorphic germ on a singular space defined by L\^e D.T. \cite{Le-oslo}. 


\begin{definition}\label{d:local}
Let $f : ( \cY, y) \dashrightarrow \bP^1$ be a germ of a meromorphic 
function. For every $a\in \bP^1$, one associates the germ $\pi_{(y,a)} : (\bY , (y, a)) \to \bP^1$. Abusing language, by restricting this map to $X\subset \bY$, we have a germ  $f: (X , (y, 
a)) \to \bC$, where the point 
$(y, a)$ might be in the closure of the set $X$.
\end{definition}

When the point $y$ does not belong to the axis $\cA$, then we have the classical situation of a holomorphic germ; the point $(y,a)$ is uniquely determined by $y$. However, when $x\in \cA$, then for each $a\in \bP^1$ we get a different germ. Proposition \ref{p:1} together with Thom's Second Isotopy Lemma show that, in this family, all germs are isomorphic except of finitely many of them.

There is a well defined local fibration at $(y,a)$, as follows.
Let $\cW$ be a Whitney stratification  of $\bY$ such that $\bY\setminus X$
is union of strata. For all small enough radii $\e$, the sphere $S_\e = 
\partial \bar B_\e(y,a)$ centered at $(y,a)$ intersects transversally all the (finitely many strata) in some neighbourhood of $(y,a)$. By \cite{Le-oslo}, the projection $\pi : \bY_{D} \cap B_\e(y,a)  \to D$ is stratified locally trivial over $D^*$, if the radius of $D$ is small enough.
It follows that the restriction:
\begin{equation}\label{eq:2}
\pi : X_{D^*} \cap B_\e(y,a)  \to D^*.
\end{equation} 

is also locally trivial. If $y$ is fixed, this fibration varies with the parametre $a$; the radius $\e$ of the ball depends also on the point $a$.
From Proposition \ref{p:1} and Thom's Isotopy Lemma it follows that, since $\pi$ is 
stratified-transversal to $\bY$ over $\bP^1 \setminus \Lambda$, the 
fibration  
$\pi : X_D \cap B_\e (y,a) \to D$ is trivial, for all but a finite number of values of $a\in \bP^1$.

\begin{definition}
We call the locally trivial fibration (\ref{eq:2})  the 
{\em Milnor-L\^e fibration}  of the meromorphic function germ $f$ {\em at the 
point $(y, a) \in \bX$}.
\end{definition}


\section{Isolated singularities and their vanishing cycles}\label{count} 

We show that, if the singularities along the indeterminacy locus
are isolated, then
one can localize the variation of topology of fibres.
The same type of phenomenon exists in the previously known cases:
 holomorphic germs \cite{Mi} and of polynomial functions (e.g. \cite[4.3]{Ti-t}).
This has consequences on the problem of detecting variation of topology, especially when the underlying space $Y$ has maximal {\em rectified homotopical depth}\footnote{This notion was introduced by Hamm and L\^e \cite{HL}
in order to realize Grothendieck's predictions that homotopical depth was the cornerstone for the Lefschetz type theorems on singular spaces \cite{Gr}. We shall come back to Lefschetz theory in \S \ref{s:lef}.}. Before stating the localization result, let us give the definition here, for further use:  Let $Z$ be a complex space endowed with some Whitney stratification $\cW$; denote by $\cW_i$ the union of strata of dimension $\le i$. After \cite{HL}, one says that $\r_\cW Z \ge m$ if for any $i$ and any point $x\in\cW_i\setminus \cW_{i-1}$, the homotopy groups of $(U_\alpha, U_\alpha \setminus \cW_i)$ are trivial up to the order $m-1-i$, where  $\{U_\alpha\}$ is some fundamental system of neighbourhoods of $x$. It is shown in {\em loc.cit.} that this doest not depend on the chosen Whitney stratification. A similar definition holds in homology instead of homotopy, giving rise to the rectified homological depth, denoted $\rHd$. Let us mention  that $\r Y \ge n$ when $Y$ is  locally a complete intersection of dimension $n$ at all its points (see \cite{LM}).
 
\begin{proposition}\label{t:loc} 
Let $F$ have isolated singularities with respect to some $\partial \tau$-stratifi- \\
cation $\cG$ at $a\in 
\bP^1$ (i.e. $\dim \bY_a \cap \Sing_\cG F \le 0$).
Then the variation of topology of the fibres of $F$ at $X_a$ is
localizable at those points.
\fin
\end{proposition}
What we mean by  ``localizable" is that there exist small enough balls in $\bY$ centered at the isolated singularities such that, outside these balls, the projection $\pi$ is trivial over a small enough disc centered at $a\in \bP^1$. The proof in the general case of meromorphic functions has the same structure as the proof presented  in \cite[Theorem 4.3]{Ti-t}, inspite the fact that in {\em loc.cit.} we consider a particular situation; we may safely skip it. 
The localization result implies, as in the particular cases \cite{Ti-t}, \cite{ST-x}, that the vanishing cycles are 
concentrated at the isolated singularities.

\begin{corollary}\label{c:conc}
Let $F$ have isolated singularities with respect to $\cG$ at $a\in \bP^1$ and let $\bY_a \cap \Sing_\cG F =\{ p_1, \ldots, p_k\}$. For any small enough balls $B_i\subset \bY$ centered 
at $p_i$, and for small enough closed disc centered at $a$, $s\in 
\partial 
D$,  we have:
\begin{equation}\label{eq:local}
 H_*(X_D, X_s)\simeq \oplus_{i=1}^k 
H_*(X_D\cap B_i, X_s\cap B_i) .
\end{equation}
\fin
\end{corollary}
 One show, following \cite{ST-x}, that an isolated $\cG$-singularity at a point of 
$A\times \bP^1\subset \bY$ is detectable by the presence of a certain local polar locus. If the space $Y$ is ``nice enough", then the local vanishing homology (second term of the isomorphism (\ref{eq:local})) is concentrated in one dimension only. Then the polar locus defines a numerical invariant which measures the number of ``vanishing cycles at this point". 
\begin{definition}\label{d:polar}
 Let $\xi= (y,a)\in A\times \bP^1$ and let $f = p/q$ a local
 representation of $F$ at $y$. Then
the {\em polar locus} $\Gamma_\xi(\pi,q)$ is the germ at $\xi$ of the 
space: 
 \[ \cl \{ \Sing_\cG(\pi, q) \setminus (\Sing_\cG \pi \cup A\times \bP^1)\} \subset 
\bY .\]
 \end{definition}

As in \cite[\S 4]{ST-x}, we have the 
isomorphisms:
\[ \Gamma_\xi(\pi,q) \simeq \Gamma_\xi(f,q) \simeq \Gamma_\xi(p,q).\]
 
 The multiplication by a unit $u$ may change the polar locus:
$\Gamma_\xi(\pi,qu)$ is in general different from $\Gamma_\xi(\pi,q)$. 
Nevertheless, 
we have the following.
\begin{proposition}\label{p:inv}
Let $\xi= (y,a)\in A\times \bP^1$ be an isolated $\cG$-singularity
of $F$ and let $f_y = p/q$ a local representation of $F$. 
Then, for any multiplicative unit $u$, the polar locus 
$\Gamma_\xi(\pi,qu)$ is either void or $\dim \Gamma_\xi(\pi,qu) =1$.
If moreover $Y$ is of pure dimension $m$ and $\r (Y \setminus \{ q=0\}) \ge m$ in the neighbourhood of $\xi$, then  the intersection multiplicity $\inter_\xi 
(\Gamma_\xi(\pi,qu), \bY_a )$ is independent on the  
unit $u$. We call this multiplicity the {\em polar number} at $\xi$. 
\end{proposition}
\begin{Proof}
We follow essentially the proof of \cite[4.2]{ST-x}. The key argument to use is the independence of $\bP T^*_{qu}$ from the multiplicative unit 
$u$, proven in \cite[Prop. 3.2]{Ti-t}, where $\bP T^*_{q}$ denotes the 
projectivized relative conormal of $q$.

Since $\dim \bP T^*_q = m+1$, it follows that $\Gamma_\xi(\pi, q)$ is either void or of dimension at least $1$. On the other hand, since $\xi$ is a point belonging to $\Gamma_\xi(\pi, q)$, it follows that 
$\Gamma_\xi(\pi, q)$ has dimension at most $1$. The same argumentation holds for $qu$ instead of $q$. This proves the first claim.

To prove the second statement, let's suppose that $\Gamma_\xi(\pi, q)$ has dimension $1$. Consider the Milnor-L\^e fibration 
 of the function $\pi$ at $\xi$:
\begin{equation}\label{eq:le}
\pi_| :  \bY_{D^*}\cap B_\e (\xi)\to D^* ,
\end{equation} 
as explained at the end of \S \ref{s:van}.
We compute the homology $H_*(\bY_s \cap B)$ of the Milnor fibre of 
this fibration. Inside $B$, the restriction  of the 
function $q$ to $\bY_s\cap B$ has a finite number of stratified isolated singularities, which are precisely the points of intersection $\bY_s \cap B \cap \Gamma(\pi,q)$.  By the result due to Goresky-MacPherson that cylindrical neighbourhoods are conical \cite[pag. 165]{GM},
it follows that $\bY_s \cap B$ 
is homotopy equivalent to $\bY_s \cap B\cap q^{-1}(\delta)$, where $\delta$ is a small disc at $0\in \bC$ such that $\bY_s \cap B \cap \Gamma(\pi, q) = q^{-1}(\delta)\cap 
\bY_s  \cap B \cap \Gamma(\pi, q)$.

 Let now take a small enough disc $\hat\delta$ centered at $0\in \bC$ such that $q^{-1}(\hat \delta)\cap 
\bY_s  \cap B  \cap \Gamma(\pi, q) =\emptyset$. By using the properties of the partial Thom stratification $\cG$ and by retraction, it follows that  $q^{-1}(\hat \delta) \cap 
\bY_s  \cap B$ is homotopy equivalent to the central fibre  $q^{-1}(0) \cap \bY_s \cap B$.  

Since the subspace $q^{-1}(0)$ at $\xi$ is the product $A\times \bP^1$,
the central fibre   
$q^{-1}(0) \cap \bY_s \cap B$ is contractible; hence $q^{-1}(\hat \delta) \cap \bY_s  \cap B$ is contractible too.

The total space $\bY_s \cap B \stackrel{\h}{\simeq} q^{-1}(\delta)\cap\bY_s \cap B$ is built by attaching 
to the space $q^{-1}(\hat \delta)\cap 
\bY_s  \cap B$  finitely many cells, which come from the isolated singularities of the function $q$ on $q^{-1} 
(\delta \setminus \hat\delta)\cap \bY_s \cap B$. Since $\r (Y_s \setminus \{ q=0\}) \ge m-1$ in some neighbourhood of $\xi$ (by our assumption and Hamm-L\^e's \cite[Theorem 3.2.1]{HL}), it follows that each singularity contributes with cells of dimension exactly $m-1$ and the number of cells is equal to the corresponding local Milnor number of 
the function $q$ (see \cite{Le-s} and \cite{Ti-b} for more details).
 The sum of these numbers is, by definition, the intersection multiplicity $\inter_\xi(\Gamma(\pi, q) , \bY_a)$.

By this we have proven that 
 $\dim H_{m-1}(\bY_s \cap B) = \inter ( \Gamma_\xi(\pi , q), \bY_a)$ and that $\tilde H_i(\bY_s \cap B) = 0$ for $i\not= n-1$. 
 Replacing $q$ by $qu$ in our proof yields the same equalities; this proves our statement. 
\end{Proof} 
\begin{remark}\label{r:spheres}
The last part of the above proof shows in fact more, that the fibre $\bY_s \cap B$ of the local fibration (\ref{eq:le})  is homotopy 
equivalent to a ball to 
which one attaches a certain number of $(m-1)$-cells. Therefore $\bY_s \cap B$ is homotopy equivalent to a bouquet of spheres $\bigvee S^{m-1}$ of dimension $m-1$.
\end{remark}

One can get more precise results when lowering the generality. The situations we consider in the following are more general that ``polynomial functions'', which we consider in \S \ref{s:poly}. Let then assume:
\vspace*{2mm}

\noindent $(*)$  $X := Y\setminus A$ has at most isolated singularities.

\vspace*{2mm}
In this case we have $\Sing_\cG F \cap A\times \bP^1 \subset \Sing \bY \cap A\times \bP^1$. In the notations of Cor. \ref{c:conc}, 
the following duality result holds (integral coefficients):
\[ 
H_*(X_D\cap B_i, X_s\cap B_i) \simeq H^{2m -*}(\bY_D \cap B_i , \bY_s \cap B_i),\]
where $m= \dim_{p_i}\bY$. This follows from Lefschetz duality for polyhedra, see Dold \cite[p. 296]{Do}.

Since $\bY_D \cap B_i$ is contractible, we get:
\begin{equation}\label{eq:dual}
 H_*(X_D\cap B_i, X_s\cap B_i) \simeq \tilde H^{2m -1 -*}(\bY_s \cap B_i) .
\end{equation}
When comparing this to Remark \ref{r:spheres} and to Definition \ref{d:vanhom}, Corollary \ref{c:conc} and Proposition  \ref{p:inv}, the following statement comes out:
\begin{theorem}\cite{ST-x}
Let $F$ have isolated singularities with respect to $\cG$ at $a\in\bP^1$ and let $\xi \in A\times \{ a\} \cap \Sing_\cG F$. Then $F_{|Y\setminus A}$ has vanishing cycles at $\xi$ if and only if $\inter (\Gamma_\xi(\pi, q), \bY_a) \not= 0$.
The number of vanishing cycles at $\xi$ is $\lambda_\xi := \dim H_{m-1} ( \bY_s \cap B) = \inter (\Gamma_\xi(\pi, q), \bY_a)
$, where $m = \dim_\xi \bY_a$.
\fin
\end{theorem}

 Let us assume now that $\bY_a$ has an isolated singularity at $\xi \in A\times \{ a\}$. This happens for instance when $Y$ has at most isolated singularities on $A$. Since $(\bY_a, \xi)$ is a hypersurface germ in $Y$, it has a well defined Milnor fibration, and in particular a Milnor number, which we denote by  $\mu(a)$.
 
 It also follows that $\dim_\xi \Sing_\cG F \le 0$ and that $\bY$ has 
singularities of dimension at most 1 at $\xi$. If $\Sing \bY$ is a curve at $\xi$, this curve intersects $\bY_s$, for $s$ close enough to $a$, at some points $\xi_i(s)$, $1\le i\le k$. There is a well defined Milnor $\mu_i(s)$ at each hypersurface germ 
$(\bY_s, \xi_i(s))$. (In case $\Sing \bY$ is just the point $\xi$, we consider that $\mu_i(s)=0$, $\forall i$.)
We have the following computation of the number of vanishing cycles at $\xi$:

\begin{theorem}\label{t:numbers}\cite{ST-x}
Let $\dim_\xi \Sing \bY_a =0$. Then:
 \[ \lambda_\xi = \mu(a) - \sum_{i=1}^{k} \mu_i(s) .\]
\end{theorem}
\begin{Proof}
 We only sketch the proof and send to \cite{ST-x} for details. Consider the function:
 \[ G = p-sq : (Y\times \bC, \xi_i(s)) \to (\bC,0).\]
 For fixed $s$, this function is a smoothing of the germ $(\bY_s,\xi_i(s))$. Since $\xi$ is an isolated singularity, the polar locus at $\xi$ of the map $(G,\pi) : Y\times \bC \to \bC^2$, defined as $\Gamma_\xi(G,\pi) = \cl \{ \Sing (G,\pi) \setminus \Sing G\}$,
is a curve. By using polydisc neighbourhoods\footnote{ polydisc neighbourhoods were first used by L\^e D.T.\cite{Le-oslo}.} $(P_\alpha \times D_\alpha)$ at $\xi$ in $Y\times \bC$, one may show that $(G,\pi)^{-1}(\eta,s) \cap (P_\alpha \times D_\alpha)$ is homotopy 
equivalent to the Milnor fibre of the germ $(\bY_a, \xi)$. To obtain from this the space $\pi^{-1}(s)\cap (P_\alpha \times D_\alpha)$, one has to attach a number of $m$ cells, where $m$ denotes $\dim_\xi Y$. Part of these cells come from the singular points $\xi_i(s) \in \Sing G \cap \pi^{-1}(s)$: by definition, their total number is $\sum_{i=1}^{k} \mu_i(s)$. The other part of the cells come from the intersection with $\Gamma_\xi(G,\pi)$ and their number is $r= \inter(\Gamma_\xi(G,\pi), \pi{-1}(0))$.
The key observation is that $r$ turns out to be equal to $\dim H_{m-1}(\bY_s \cap B)$, which is the number of vanishing cycles at $\xi$.
Finally, since $\pi^{-1}(s)\cap (P_\alpha \times D_\alpha)$ is contractible (since being the Milnor fibre 
of a linear function $t$ on a smooth space), we have the following equality:

 \[  \mu(a) = r + \sum_{i=1}^k \mu_i(s).\]
Notice that, in case $\dim_\xi \Sing \bY =0$, we get just $\lambda_\xi = \mu(a)$.
\end{Proof}
We send to Corollary \ref{c:num} for the counting of the total number of vanishing cycles. 
Let us give two examples of meromorphic functions, one on $\bP^2(\bC)$ and another on a nonsingular quadratic surface 
in  $\bP^3(\bC)$. 

\begin{example}\label{ex:1}(\cite{ST-x})
Let $F \colon Y=\bP^2 \dashrightarrow \bP^1$,   $\displaystyle F= \frac{x(z^{a+b} + 
x^ay^b)}{y^pz^q}$, where $a+b+1 = p+q$ and $a,b,p,q \ge 1$. 
 For some $s\in \bC = \bP^1\setminus [1:0]$, the space $\bY_s$ is given by:
 \begin{equation}\label{eq:1}
 x(z^{a+b} + x^ay^b) = s y^pz^q
 \end{equation}
 $\Sing \bY \cap (Y\times \bC)$ consists of three lines: $\{ [1:0:0], [0:1:0], [0:0:1]\} 
\times \bC$. We are under the assumptions of Theorem \ref{t:numbers} and we inspect each of these 3 families of germs with isolated singularity to see where the Milnor number jumps.
 
 Along $[1:0:0]\times \bC$, in chart $x=1$, there are no jumps, since the germs have uniform Brieskorn type $(b, a+b)$.
 Along $[0:0:1]$, in chart $z=1$, there are no jumps, since the type is constant $A_0$, for all $s$.
 Along $[0:1:0]$, in chart $y=1$: For $s\not= 0$, the Brieskorn type is $(a+1, 
q)$, with $\mu(s) = a(q-1)$. If $s=0$, then we have $x^{a+1} + xz^{a+b}=0$ 
with $\mu(0) = a^2 + ab + b$.

 There is only one jump, at $\xi = ([0:1:0], 0)$; according to the preceding theorem,
$\lambda_\xi = a^2 + ab + b - a(q-1) = b +ap$.
 
\end{example}



\section{Homotopy type of fibres}


In case of 
isolated singularities, we have the following result on the relative homotopy type.
\begin{theorem}\label{t:attach} 
 Let  $\{ a_i\}_{i=1}^p$ be the set of atypical values of $F$ within some  
 open disc $D\subset \bP^1$. Let $s\in D$ be some typical value of $F$.
 For all $1\le i\le p$, let  $F$ have an isolated $\cG$-singularity at $a_i$,  $Y$ be of pure dimension $m$ at $a_i$ and $\rHd (Y\setminus (V\cup A)) \ge m$ in some neighbourhood of $a_i$. If either of the two following conditons
 is fulfilled:
\begin{enumerate}
 \item $X$ is a Stein space,
 \item $X= Y\setminus A$ and $X$ has at most isolated singularities,
\end{enumerate} 
  then $X_D$
 is obtained from 
$X_s$ by attaching cells of real dimension $m$.
In particular, the topological space $X_D/X_s$ is homotopy equivalent to a bouquet of speres $\vee S^m$.
\end{theorem}
\begin{Proof} 
We prove that the reduced integral homology of $X_D/X_s$ is concentrated in dimension $m$. 
By Proposition \ref{t:loc}, the 
variation of topology of the fibres of $F_{| X_D}$ is
localizable at the points $\Sing_\cG F \cap \bY_D$. 
We have to take into account all the possible positions of such a singular point $\xi= (y,a)$, namely: on $X$, on $V\setminus A$ or on $A\times \bP^1 \subset \bY$. 

In all the cases, it turns out that the pair $(X_{D_a}\cap B_\xi, X_s \cap B_\xi)$ is ($m -1$)-connected, where $B_\xi \subset \bY$ is a small enough ball at $\xi$, $D_a$ is a small enough closed disc at $a$ and $s\in \partial D_a$.

For a point $\xi$ in the first case, this is just Milnor's classical 
result for holomorphic functions with isolated singularity 
\cite{Mi}.
In the two remaining cases, this follows by a result due to
 Hamm and L\^e \cite[Corollary 4.2.2]{HL}, in a slightly improved version for partial Thom 
stratifications (see \cite[2.7]{Ti-t}). This result needs the condition on the rectified homological depth.

  By the above proven connectivity  
 of the pair $(X_{D_a}\cap B_\xi, X_s \cap B_\xi)$ and the splitting of 
 vanishing homology into local contributions Corollary \ref{c:conc}, we get that the homology of $(X_D,X_s)$ is zero below dimension $m$.
 Above dimension $m$, we also have the annulation, due to the following reasons. In case (a) the space $X_D$, respectively $X_s$, is Stein of dimension $m$, resp. $m-1$. In case (b),  we may apply the duality (\ref{eq:dual}) and we have that the cohomology 
 $\tilde H^*(\bY_s\cap B_\xi)$ is concentrated in dimension $m-1$, by Remark \ref{r:spheres}.
 
  Then one can map a bouquet of $m$ spheres into $X_D/X_s$
such that this map is an isomorphism in homology. This implies, by Whitehead's theorem (see \cite[7.5.9]{Sp}), that the map induces an isomorphism of homotopy groups. (Remark that $X_D/X_s$ is simply connected whenever $m\ge 2$). Since we work with analytic objects, therefore triangulable, the space $X_D/X_s$ is a CW-complex. We may now conclude our proof, since for CW-complexes, weak homotopy equivalence coincides with homotopy equivalence. 
\end{Proof}

In case (b), this result has been proved by Siersma and the author \cite{ST}. Let us point out that, in this case (b), we also have the local bouquet result: $(X_{D_a}\cap B_\xi) /  (X_s \cap B_\xi) \stackrel{\h}{\simeq} \bigvee S^m$. 

When assuming high connectivity of the space $X$, we get the following immediate consequence (proved in lower generality in \cite{ST-x}).
\begin{corollary}
Under the hypotheses of Theorem \ref{t:attach},  if in addition the space $X$ is Stein and $(m-1)$-connected, 
then $X_s \stackrel{\h}{\simeq} \bigvee S^{m-1}$.
\fin
\end{corollary}
  Particular cases of this corollary appeared previously in several circumstances: Milnor's bouquet result \cite{Mi} on holomorphic germs with isolated singularity;
bouquet results for generic fibres of polynomial maps with isolated singularities in the affine \cite{Br-1,Br} and with isolated singularities at infinity \cite{ST}, \cite{Ti-t}.

 As another consequence, we shall draw a formula for the total number of
 vanishing cycles in case of isolated 
$\cG$-singularities. Let us denote by $\lambda_a$ the sum of the polar  numbers at the singularities on $(A\times \bP) \cap \bY_a$ and by 
$\mu_a$ the sum of the Milnor numbers of the singularities on $\bY_a \setminus (A \times \{a\})$. One needs to note that the Milnor fiber (in our case $\bY_s \cap B$) of a holomorphic function with isolated singularity on a Whitney stratified space is homotopically a bouquet of spheres of dimension $=\dim \bY_s \cap B$, see \cite{Ti-b}.
\begin{corollary}\label{c:num}
 Under the hypotheses of Theorem \ref{t:attach}, we have:
\[  \dim H_{m-1}(X_{D_a}, X_s) =  \mu_a + \lambda_a \ , 
\]
\[ \dim H_{m-1}(X_D, X_s) = \sum_{a\in D}\mu_a + \sum_{a\in D}\lambda_a . \]
\fin
\end{corollary}
\begin{remark}\label{r:ex} 
In Example \ref{ex:1}, let us consider $X=\bP^2\setminus \{ yz=0\}$. Let $s\in \bC\subset \bP$.
 It is easy to see that the fiber $X_0$ is a disjoint union of $c+1$ 
disjoint copies of $\bC^*$, where $c= \gcd (a,b)$, therefore $\chi(X_0) =0$. For $s\not= 0$, by a branched covering 
argument, one shows $\chi(X_s) =- (b+ap)$. The vanishing homology is concentrated in dimension 2, by Theorem \ref{t:attach}. 
 When taking $D=\bC$, we get the Betti number
  $b_2(X,X_s) = \chi(X,X_s) = \chi(X) - \chi(X_s) = (b+ap)$.
  We have seen at \ref{ex:1} that the sum $\sum_{a\in \bC}\lambda_a$ 
  consists of a single term $\lambda_\xi = b+ap$. One can also see easily that there is no other singularity, in particular that $\sum_{a\in \bC}\mu_a =0$. Hence the second equality in Corollary 
\ref{c:num} is verified.
  
\end{remark}

\section{Monodromy}
  

Let $F\colon Y \dashrightarrow \bP^1$ be a meromorphic function and let $\Lambda\subset \bP^1$ denote the set of atypical values of the associated map $\pi \colon \bY \to \bP^1$. 
There is a well defined monodromy $h_i$ around an atypical value $a_i \in 
\Lambda$.  This is induced by a counterclockwise loop around the small circle $\partial D_i$. Let $D$ denote some large disc, like in \S \ref{s:van}, such that $D\cap \Lambda \not= \emptyset$. We have a geometric monodromy representation:
\[ \rho_i : \pi_1 (\partial \bar D_i, s_i) \to \Iso (X_D, X_{\partial \bar D_i},X_{s_i}), \]
where $\Iso (.,.,.)$ denotes the group of relative isotopy classes of stratified homeomorphisms (which are C$^\ity$ along each stratum).
Note that the retriction of this action to $X_D$, $X_{D_i}$ or to $X_{\partial D_i}$,
is trivial. Let $T_i$ denote the action induced  by $\rho_i$ in homology
(with integral coefficients).

Let us identify $H_*(X_D,X_s)$ to 
$\oplus_{a_i\in D\cap \Lambda} H_*(X_{D_i}, X_{s_i})$ as in Proposition \ref{p:basic}. This identification depends on the chosen system of paths $\gamma_i \subset D$ from $s\in \partial D$ to $s_i\in \partial D_i$, as explained in \S \ref{s:van}. 
There is the following general result, showing that the 
action of the monodromy $T_i$ on a vanishing cycle $\omega \in H_*(X_D,X_s)$ 
changes $\omega$ by adding to it only contributions from the homology 
vanishing at $a_i$. 

\begin{proposition}\label{p:picard}
For every $\omega\in H_*(X_D,X_s)$, there is $\psi_i(\omega) \in H_*(X_{D_i}, X_{s_i})$ such that 
$T_i (\omega) = \omega + \psi_i(\omega)$.
\end{proposition}
\begin{Proof}
The proof goes exactly as in the more particular that we consider in \cite[Prop. 6.1]{ST-x}. 
One may identify the map:
$T_i - \id : H_{q+1}(X_D,X_s) \to H_{q+1}(X_D,X_s)$
to the composed map:
\begin{equation}\label{eq:wang} 
H_{q+1}(X_D,X_s) \stackrel{\partial}{\to} H_q(X_s) \stackrel{w}{\to}
H_{q+1}(X_{\partial D_i},X_s) \stackrel{i_*}{\to} H_{q+1}(X_D,X_s),
\end{equation}
where $w$ denotes the Wang map, which $w$ is an isomorphism, by K\"unneth 
formula. The last morphism in (\ref{eq:wang}) factors as follows:
\[ \begin{array}{rcl}
  H_{q+1}(X_{\partial D_i},X_{s_i}) & \stackrel{i_*}{\longrightarrow}  & 
H_{q+1}(X_D,X_s) \\
  \searrow & \  & \nearrow  \\
 \ &  H_{q+1}(X_{D_i}, X_s) & \ 
\end{array} \] 
where all three arrows are induced by inclusion. It follows that 
the submodule of ``anti-invariant cycles'' $\im (T_i - \id : 
H_* (X_D,X_s)\to H_* (X_D,X_s))$ is 
contained in the direct summand $H_*(X_{D_i},X_{s_i})$ of $H_*(X_D,X_s))$.
\end{Proof}
One has the following easy consequence, in full generality. Assume that the paths in $D$, say $\gamma_1, \ldots 
\gamma_l$ are counterclockwise ordered. The chosen paths define a  decomposition of $H_*(X_D,X_s)$ into the direct sum $\oplus_{a_i\in D\cap \Lambda} H_*(X_{D_i}, X_{s_i})$. Denote by $T_{\partial D}$ the monodromy around the circle $\partial D$. One has the following imediate consequence,  remarked in \cite{DN2} for polynomial functions and in \cite{ST-x} for the particular case $Y$ nonsingular and $X=Y\setminus A$.
\begin{corollary}\label{c:pis}
 Assume that the direct sum decomposition of $H_*(X_D,X_s)$ is fixed. Then $T_{\partial D}$ determines $T_i$, $\forall 
i\in\{ 1, \ldots , l\}$.
\fin
\end{corollary}

\begin{note}\label{n:piclef}
One may say that Proposition \ref{p:picard} is a Picard type formula, since Picard showed it at the end of the XIX-th century, for algebraic functions of two variables with simple singularities. Lefschetz proved later the well known formula for a loop around a quadratic singularity, in which case $\psi_i(\omega)$ is, up to sign, equal to $c\Delta$, where $\Delta$
is the quadratic vanishing cycle and $c$ is the intersection number $(\omega, \Delta)$. This became the basis of what one calls now Picard-Lefschetz theory (which is the conterpart of the Morse theory, in case of complex spaces), see e.g. \cite{AGV}, \cite{Eb}, \cite{Va}. In case of polynomial functions, the Picard formula was singled out in \cite{DN1}, \cite{NN2}.
\end{note}   
\begin{remarks}\label{r:piclef}
 The statement and proof of Proposition \ref{p:picard} dualize easily from homology to cohomology. One 
obtains in this way statements about invariant cocycles $\ker (T^i -\id 
\colon H^*(X,F) \to H^*(X,F))$ instead of anti-invariant cycles.

  A special case is that of a polynomial function 
$F : \bC^n \to \bC$, for which $X = \bC^n$. 
Results on invariant cocycles were obtained in 
\cite{NN}; they can be proved also in our more general setting. 

\end{remarks}

In the rest of this section we review some results on the zeta function
of the monodromy. We shall only discuss  global meromorphic functions $F$; following the general remark in the Introduction, all results translate easily in case of meromorphic germs.
For the particular case of polynomial functions, we send the reader
to \S \ref{s:poly}, where we present more specific results. 

\begin{definition}
 Let $T_{\partial D}$ be the monodromy around some disk $D$ as above. One calls zeta function of $T_{\partial D}$ the following rational function in variable $t$:
\[ \zeta_{(X_D, X_s)} (t) = \prod_{i\ge 0} \det [\id - t T_{\partial D} : H_i(X_{D},X_s)\to H_i(X_{D},X_s)
]^{(-1)^{i+1}}.\]
\end{definition}

 We are interested here in the zeta function  of the monodromy around a value
 $a\in \Lambda$. Let us first assume that $F$ has isolated
 $\cG$-singularities. By the direct sum splitting (Corollary \ref{c:conc})
 and since the monodromy acts on each local Milnor fibration, we get:
 \[ \zeta_{(X_{D_a}, X_s)} (t) = \prod_{i=1}^k \zeta_{(X_{D_a} \cap B_i, X_s \cap B_i)} (t),\]
 where $\{ p_1, \ldots, p_k\} = \bY_a\cap \Sing_\cG F$ and $B_i$ is a small Milnor ball centered at $p_i$.
 
 For the zeta function  of the monodromy $T_{\partial D_a}$ acting on the homology of the general fibre $X_s$ we also get:
\[ \zeta_{X_s} (t)= (1-t)^{-\chi(X_a)} \prod_{i=1}^k \zeta^{-1}_{(X_{D_a} \cap B_i, X_s \cap B_i)} (t),\]
 since the monodromy acts on $X_{D_a}$ as the identity and since $\chi(X_{D_a}) = \chi(X_a)$.

 Let us now suppose that $Y$ is nonsingular and consider $X= Y\setminus A$,
 but not assume anything about the singularities of the meromorphic function $F$.
 One may follow the method of A'Campo \cite{A'C-zeta} to produce a formula for  the zeta function, as follows. There exists a proper holomorphic modification
$\phi \colon \tilde Y \to Y$,  which is bi-holomorphic over $X\setminus \cup_{a\in \lambda} X_a$.
 The pull-back $\tilde F = F\circ \phi$ has a general fibre $\tilde X_s$ which is isomorphic to $F$. The action of the monodromy is also the same, therefore $\zeta_{\tilde X_s} (t) = \zeta_{X_s} (t)$. Then one can write down a formula for the zeta function around the value $a\in \bP^1$ in terms of the exceptional divisor and the axis $A$. By expressing the result as an integral with respect to the Euler characteristic (see e.g. Viro's paper \cite{Vi} for this technique), one can get rid of the resolution.
\begin{proposition}\label{p:glm} \cite{GLM-1} 
 Let $X= Y\setminus A$ and let $Y$ be nonsingular. Then: 
 \[ \zeta_{X_s}(t) = \int_{A\times \{ a\}\cup X_a} \zeta_p (t) d\chi,\]
 where $\zeta_p$ denotes the local zeta function at the point $p\in \bY$.
\fin 
\end{proposition}
 Further formulae for the zeta function and some consequences can be found in the papers by Gusein-Zade, Luengo and  Melle \cite{GLM-1, GLM, GLM-3}.
 
  From $\zeta_{X_s}$ one easily gets $\zeta_{(X_D, X_s)}$, since $\zeta_{X_s}= \zeta_{X_D}\cdot \zeta_{(X_D, X_s)}^{-1}$ and $\zeta_{X_D} = (1-t)^{-\chi(X_a)}$.



\section{Nongeneric pencils and Zariski-Lefschetz type results}\label{s:lef}

Exploring a space $Y$ by pencils of complex hyperplanes is an old idea 
in mathematics. It appeared in Lefschetz's work \cite{Lef}, which became the fundation of the so-called Lefschetz Theory. Almost in the same time the Morse Theory was born \cite{Mo}.
Each of the two theories provide a method for studying the topology of the space; both use scanning with levels of a function.\footnote{for hystorical notes and new developments until about 1987, see Goresky and MacPherson's  book \cite{GM}} 
The analogous of Morse function for the  Lefschetz Theory is ``Lefschetz pencil''.

One usually means by Lefschetz pencil a pencil having
singularities of simplest type (i.e.  A$_1$) and transversal axis. In the usual projective space, such pencils are generic, but on certain spaces they might not even  exist. A more general point of view is to allow pencils with isolated singularities, alias meromorphic functions $F \colon Y \dashrightarrow \bP^1$ with isolated singularities in the sense of this paper. We call them ``nongeneric pencils" and point out that they can have singularities also within the axis $A$ of the pencil. The case of isolated singularities outside the axis has been considered before by Hamm and L\^e (e.g. \cite{HL}) and by Goresky and MacPherson (see \cite{GM}).

 The following connectivity result of Zariski-Lefschetz type holds. 
\begin{theorem}\label{t:lef} \cite{Ti-lef}
Let the pencil $F= P/Q \colon Y \dashrightarrow \bP^1$ have isolated $\cG$-singularities (Definition \ref{d:sing}).
Assume that 
$A\not\subset V$ and let $X_s$ denote a generic member of the pencil.

If $\r X \ge m$, $m\ge 2$, and if the pair $(X_s , A\cap X_s)$ is 
$(m-2)$-connected 
then the morphism induced by inclusion:
\[ \pi_i (X_s) \to \pi_i (X) \]
is an isomorphism for $i\le m-2$ and an epimorphism for $i=m-1$.
 \fin
\end{theorem} 

This represents a far-reaching extension of the classical Lefschetz 
theorem on hyperplane sections. The latter says that, if $X$ is a projective variety and $H$ is a hyperplane such that $X\setminus H$ is nonsingular of dimension $m$, then $\pi_i(X,X\cap H)= 0$, for all $i\le m-1$. This can be viewed as a statement about pencils with transversal axis (i.e. there are no singularities on $A$), since, even if $H$ is not a generic hyperplane, it is a member of some generic pencil in the projective space. Indeed, one may define 
such a pencil by  
choosing a generic axis inside the hyperplane $H$. Then our claim follows by the conjunction of the following 2 observations: 1). Theorem \ref{t:lef} is true when replacing $X_s$ by the tube $X_{D_a}$, where one  allows singularities of any type on $Y_a$; 2). $X_{D_a}$ is contractible to $X_a$ when $X$ is compact (i.e. $V=\emptyset$).

 One may draw the following consequence on the homotopy type of the pair space-section, which actually represents an extension of Theorem \ref{t:attach}:

\begin{corollary}\label{c:lef} \cite{Ti-lef}
 Under the hypotheses of Theorem \ref{t:lef}, up to homotopy type, 
 the space $X$ is built 
from $X_s$ by attaching cells of dimension $\ge m$. 
If $X$ 
is in addition a Stein space of dimension $m$, then the attaching 
cells are of dimension precisely $m$.
\fin
\end{corollary}

\begin{remark}\label{r:lef}
What happens when $A\subset V$? 
We prove in \cite{Ti-lef} that if $\{ Q=0\}\subset V$ (which is a special case of $A\subset V$, since $A\subset\{ Q=0\}$), then the conclusion of Theorem \ref{t:lef} holds, with the single assumption $\r X \ge m$.

This result concerns in particular 
the polynomial functions $P : \bC^n \to \bC$. Such a function defines a nongeneric pencil on $X= \bC^n$, since it can be regarded as a meromorphic 
function $\tilde P/ Q$ on the (weighted) projective space $Y = \bP_w^n$, where  $\{ Q=0\}$ is the hyperplane at infinity. Notice that in this case the condition $\r X \ge n$ holds since $X$ is nonsingular. We refer to \S \ref{s:poly} for some consequences.
\end{remark}

\section{Equisingularity at the indeterminacy locus}\label{s:equi}

Equisingularity conditions are considered beginning with  Zariski's work
on families of algebraic (hyper)surfaces, see \cite{Za}. There are more recent contributions to local equisingularity theory, especially by Teissier \cite{Te-1, Te-2} and Gaffney (\cite{Ga} and several other papers of the same author).

The case of families on non-compact spaces, like our family $\{ X_a\}_{a\in \bP^1}$,  where $X = Y\setminus V$ and $A\cap V \not= \emptyset$,
is special. We have seen 
 that the singularities at the indeterminacy locus play an important role.
 The problem would be to find the weakest equisingularity condition at $A$ such that to imply topological triviality in the neighbourhood of $A$.
  If we stratify everything by Whitney conditions, then we may invoke Whitney equisingularity, which implies topological triviality; but Whitney equisingularity is too strong.  The search for a weaker alternative 
 has itself some history behind; maybe the first result in this sense is L\^e-Ramanujam's theorem for families of holomorphic germs with isolated singularity: ``$\mu$ constancy implies topological triviality'', see \cite{LR}.
  It had been found that $\mu$ constancy is really weaker than Whitney equisingularity \cite{BS} (i.e. that $\mu$ constancy is weaker than $\mu^*$ constancy\footnote{see Teissier's paper \cite{Te-2} for $\mu^*$ constancy}).
 
 In the same spirit, the problem was formulated (and solved) in case of a family of affine hypersurfaces in \cite{Ti-e}, where the equisingularity at infinity, respectively C$^\infty$-triviality at infinity, comes into the picture (see next section for details). In {\em loc.cit.}, these two notions are related to the partial Thom stratification $\cG$.
 For a family defined by a meromorphic function $F \colon Y \dashrightarrow \bP^1$, one may follow the ideas of \cite{Ti-e} up to some extent, as we pointed out in \cite{ST-x}. Let us give the main lines.  
 
\begin{definition}
 We say that $F_{|X}$ is topologically trivial at $\xi\in A\times \bP^1 \subset \bY$, resp. at $a\in \bP^1$,  if there is a neighbourhood
  $\cN$ of $\xi$, resp. of $\bY_a \cap (A\times \bP^1)$, and a small enough 
disc $D$ at $a$ such that the map $\pi_| : \cN \cap X_D \to D$ 
is a trivial fibration. 
\end{definition}

The points $\xi\in A\times \bP^1$ which pose problems are those in
$\Sing_\cG F$, since for the others we have the topological triviality.
 This claim follows by attentive re-reading of  
 \cite[Theorems 2.7, 4.6, 1.2]{Ti-e}, \cite[Theorem 7.2]{ST-x}; we actually get the following:

\begin{theorem}\label{t:equi}
Let $F = P/Q$ have isolated $\cG$-singularities at $\xi$, resp. at $a\in \bP^1$. Assume that $Y$ is of pure dimension $m$, that $X= Y\setminus \{ Q=0\}$ and that $\r X \ge m$.
 
Then $F_{|X}$ is  topologically trivial at $\xi$, resp. at 
$a\in \bP^1$ if and only if $\lambda_\xi = 0$, resp. $\lambda_a = 0$.

In particular, if $X$ has isolated singularities and $F$ has isolated $\cG$-singularities at $a\in \bP^1$, then 
$X_a$ is a general fibre of $F_{|X}$ if and only if $X_a$ is nonsingular and $\lambda_a = 0$.
 \fin
\end{theorem} 

Combining Theorem \ref{t:equi} with Corollary \ref{c:num}, we get 
the following consequence: 

\begin{corollary}\label{c:jump}
Under the hypotheses of Theorem \ref{t:attach}, a fibre $X_a$ of $F_{|X}$ is  general if and only if it has the same Euler characteristic of a general fibre.
\fin
 \end{corollary}
 
 Both results above have been stated, in
slightly lower generality, in \cite{ST-x}.
Corollary \ref{c:jump} extends the criteria for atypical 
fibres in case of polynomial functions in 2 variables \cite{HaLe}, and in 
$n$ variables \cite{ST}, \cite{Pa}. See also Proposition \ref{p:cond}. 


\section{More on polynomial functions}\label{s:poly}

In the last years there has been developed a flourishing activity in the topology of
polynomial maps, partly due to the links with affine geometry (see e.g. Kraft's Bourbaki talk \cite{Kr}). We give here a brief overview, throughout some of the multitude of the contributions.

We have explained in the Introduction that a complex polynomial function $P: \bC^n \to \bC$, $\deg P =d$, can be extended to a meromorphic function $\frac{\tilde P}{x_0^d} : \bP^n \dashrightarrow \bP$. Here $X=\bC^n$ and $Y = \bP^n$, are nonsingular spaces. This is not the only way of extending $P$ and the space $\bC^n$; one may consider\footnote{see \cite{Ti-t} for a general treatment} for instance an embedding of $\bC^n$ into some toric variety, such as a weighted projective space $\bP^n_w$.

Maybe the first author who studied the topology of polynomial functions was Broughton \cite{Br-1}. In the same time Pham \cite{Ph} was interested in regularity conditions under which a polynomial has good behaviour at infinity. 
 Some of the challenging problems that have been under research ever since are:
\vspace*{2mm}

\noindent
 1. Determine the smallest set $\Lambda$ of atypical values of $P$.\\
 2. Describe the topology of the general fibre and of the atypical fibres.\\
 3. Describe the variation of topology in the family of fibres; monodromy.\\

 In problem 1., there are only partial answers.  One has to decide which are the atypical values among a finite set of values singled out by Proposition \ref{p:1}. For instance, our general result Corollary \ref{c:jump} applies here (see also the comment following it). It is easy to show that singular values of $P$ are atypical. Then, fixing a nonsingular fibre $X_a$ of $P$, one may try to prove topological triviality at infinity by constructing a controlled vector field and ``pushing'' $X_a$ along it. This is an idea due to Thom \cite{Th}.
  There are mainly two strategies: to work in the Euclidean space or to 
  compactify $\bC^n$ into some $Y$ and use the stratification $\cG$.
  
  The first one leads to regularity conditions, in more and more generality: tameness \cite{Br}, quasi-tameness \cite{Ne}, M-tameness \cite{NZ}, $\rho$-regularity \cite{Ti-r}.
  There are also the Malgrange condition (see \cite[2.1]{Ph})---which is a condition on the \L ojasiewicz number at infinity---and its generalization by Parusi\'nski \cite{Pa}. 
  
  The second strategy leads to the t-regularity \cite{ST}, or more generally, equisingularity at infinity (which has been discussed for meromorphic functions in \S \ref{s:equi}).
    
 There are of course relations between all these conditions; one may consult \cite{Ti-r} and its references.
  Under certain circumstances, several of these conditions are equivalent.
  We may quote the following result:
\begin{proposition}\label{p:cond} \rm (after \it \cite[5.8]{ST} \rm and  \cite[2.15]{Ti-r})\it \\ 
Let $P \colon \bC^n \to \bC$ have isolated $\cG$-singularities at infinity, at $a\in \bC$.
Then the following are equivalent:
\begin{enumerate}
\item $P$ is M-tame at $X_a$.
\item $P$ satisfies Malgrange condition at infinity at $X_a$.
\item $P$ is t-regular at infinity at $X_a$.
\item $\lambda_a =0$.
\end{enumerate}
\end{proposition}
 It follows that, for a polynomial with isolated $\cG$-singularities at infinity, a fibre $X_a$ is general if and only if $\mu_a =0$ and one of the above conditions are satisfied.
 
   In case of 2 variables, the hypothesis of the above statement \ref{p:cond} is fulfilled and therefore the conclusions are valid. Indeed, 
in 2 variables, any reduced fibre $X_a$ has at most isolated $\cG$-singularities at infinity. Moreover, there are several other criteria expressing non singularity at infinity, equivalent to the ones above; we send the reader to \cite{Du}, \cite{Ti-r}.

 Still in 2 variables, one may derive the following equivalent formulation of the well known Jacobian Conjecture, in terms of singularities at infinity \cite{LW,ST}: {\em Let $P\colon \bC^2 \to \bC$ be a polynomial without critical point.
 If there exists a value $a\in \bC$ such that $\lambda_a \not= 0$, then for any other polynomial $g$, the zero locus of their jacobian ideal $\Jac(P,g)$ is not void.}

 Revisiting equisingularity at \S \ref{s:equi}, one gets more specific results in case of polynomials. By taking hyperplane sections with respect to $P$ in $\bC^n$, one may define global polar curves \cite{Ti-e} and use them in order to define an intersection number with some fibre $X_a$. By restricting $P$ to a general hyperplane and repeating generical cutting, we get a sequence of intersection numbers\footnote{they are invariant under linear change of coordinates but not under affine automorphisms.} $\gamma_a^{n-1}, \ldots , \gamma_a^{0}$. We show in \cite[Theorem 1.1]{Ti-e} that the constancy of $\gamma_s^*$, for $s$ in some neighbourhood of $a$, is equivalent to the equisingularity at infinity of $p$ at $X_a$.
 
 Secondly, the slicing processus just described gives a model of a fibre  $X_a$ as CW-complex. Let $\lambda_a^i := \gamma_s^i - \gamma_a^i$, where $s$ is a typical value of $P$.
\begin{theorem}\label{t:cw}\cite{Ti-e}
 Let $P \colon \bC^n \to \bC$ be a polynomial function. 
 Suppose that the fibre $X_a = P^{-1}(a)$ has at most 
isolated singularities. 
Then $X_a$ is homotopy equivalent to a generic hyperplane 
section $X_a \cap H$ to which one attaches $\gamma_a^{n-1} - \mu(X_a)$ cells of dimension $n-1$.

Moreover, $X_a$ is homotopy equivalent to the CW-complex obtained 
by successively attaching 
to $\gamma_a^{0} =\deg P$ points a number of $\gamma_a^1$ cells of dimension 1, then 
$\gamma_a^2$ cells of dimension 2, $\ldots$, $\gamma_a^{n-2}$ cells of 
dimension $n-2$ and finally $\gamma_a^{n-1} - \mu(X_a)$ cells of dimension 
$n-1$.
 In particular,
 $\chi(X_a) = (-1)^n \mu(X_a) + \sum_{i=0}^{n-1} (-1)^i\gamma_a^i$
and\\
 $\chi(X_s) - \chi(X_a) = (-1)^{n-1} \mu(X_a) + \sum_{i=0}^{n-1} 
(-1)^i\lambda_a^i$.
  \fin
\end{theorem} 
 One may compare this result to Corollary \ref{c:num} and notice that the
 sequence $\lambda_a^*$ is a refinement of the number of vanishing cycles at infinity $\lambda_a$, in case 
 of isolated $\cG$-singularities. Nevertheless, the numbers $\lambda_a^i$
 are defined without any hypothesis on singularities at infinity.
 
 The vanishing cycles at infinity were described for the first time in \cite{ST}. It was shown in {\em loc.cit.} that, for $P$ with isolated $\cG$-singularities at infinity, the vanishing homology is concentrated in dimension $n$; this implies that the general fibre is homotopically a bouquet of $(n-1)$-spheres.  Further progress in describing the general fibre, the special fibres and the vanishing cycles was made by Neumann-Norbury \cite{NN, NN2}, Dimca-N\'emethi \cite{DN1}, the author \cite{Ti-lef}. The cohomology of fibres is investigated by Hamm \cite{Ha-m}.

The above construction of the model for the general fibre can be 
pushed further; one may construct a global geometric monodromy group, acting on this model. This yields localization results and formulae for the zeta-function \cite{ST-m}. The more geometric point of view
on monodromy at a singularity at infinity gives two types of singularities
with local $\lambda$ equal to 1, see \cite[\S 6]{ST-m}. This may be contrasted with $\mu =1$ in case
of holomorphic germs, when the singularity can only be of one type, A$_1$. 

  One of the monodromies in case of a polynomial $P$ is the one around a big disc containing all the atypical values, denoted $T_\ity$. In two variables, 
 Dimca \cite{Di} shows that $T_\ity$ acting on the cohomology of the fibre is the identity if and only if the monodromy group of $P$ is trivial; the eigenvalue 1 occurs only in size one Jordan blocks. We send the reader to {\em loc.cit.} for further results and their discution in contrast to the holomorphic germs case. 
 
 Further aspects, such as mixed Hodge structure on fibres and algebraic Gauss-Manin systems, have been studied by several authors: Garc\'{\i}a-L\'opez and N\'emethi \cite{GN1, GN2}, respectively
 Dimca-Saito \cite{DS}, Sabbah \cite{Sa1, Sa2}.

\end{document}